\documentclass[final,leqno,onetabnum]{siamltex704}
\usepackage{amsmath} 
\usepackage{amsfonts} 
\usepackage{exscale}

\usepackage{ifpdf}
\ifpdf
\usepackage[%
  pdftitle={Bounds on changes in Ritz values for a perturbed invariant subspace
of a Hermitian matrix},%
  pdfauthor={M. E. Argentati, A. V. Knyazev, C. C. Paige, I. Panayotov},%
  pdfsubject={Rayleigh-Ritz method error bounds},%
  pdfkeywords={Hermitian matrices, angles between subspaces,
principal angles, canonical angles, majorization, Lidskii's
eigenvalue theorem, perturbation bounds, Ritz values, Rayleigh-Ritz
method, invariant subspace},%
  pdfstartview=FitH,%
  bookmarks=true,%
  bookmarksopen=true,%
  breaklinks=true,%
  colorlinks=true,%
  linkcolor=blue,anchorcolor=blue,%
  citecolor=blue,filecolor=blue,%
  menucolor=blue,pagecolor=blue,%
  urlcolor=blue]{hyperref}
\else
\usepackage[%
  breaklinks=true,%
  colorlinks=true,%
  linkcolor=blue,anchorcolor=blue,%
  citecolor=blue,filecolor=blue,%
  menucolor=blue,pagecolor=blue,%
  urlcolor=blue]{hyperref}
\fi

\usepackage{exscale}
\renewcommand{\ldots}{\ensuremath{\dotsc}}

\newcommand{\be}{\begin{equation}}
\newcommand{\ee}{\end{equation}}

\newcommand{\ba}{\begin{array}}
\newcommand{\ea}{\end{array}}

\newcommand{\da}{{\downarrow}}
\newcommand{\ua}{{\uparrow}}
\newcommand{\wsm}{{\ \prec}_{w}\ }

\newcommand{\C}{\mathbb{C}} 

\newcommand{\hm}{\mathcal{H}(n)}
\newcommand{\um}{\mathcal{U}(n)}
\newcommand{\Ra}{\mathcal{R}} 
\newcommand{\R}{\mathbb{R}} 
\newcommand{\X}{\mathcal{X}}
\newcommand{\Y}{\mathcal{Y}}
\newcommand{\spr}{{\rm{spr}}}
\newcommand{\gap}{{\rm{gap}}}
\newcommand{\ones}{{e}}

\newtheorem{conj}{Conjecture}[section]

\title{Bounds on changes in {R}itz values for a perturbed invariant subspace
of a {H}ermitian matrix\thanks{
Received by the editors July 10, 2006;
accepted for publication (in revised form) January 25, 2008;
published electronically ???????????.
}}

\author{M. E. Argentati
\footnotemark[2]\ \footnotemark[3]
\and
A. V. Knyazev\footnotemark[2]\ \footnotemark[4]
\and
C. C. Paige\footnotemark[5]
\and  I. Panayotov\footnotemark[6]
}

\begin{document}
\vspace{-1.2in}
\slugger{simax}{2008}{00}{0}{000-000} 
\vspace{.9in}

\setcounter{page}{1}
\maketitle

\renewcommand{\thefootnote}{\fnsymbol{footnote}}
\footnotetext[2]{Department of
Mathematical Sciences, University Colorado Denver, 
P.O. Box 173364, Campus Box 170, Denver, CO 80217-3364, USA.}
\footnotetext[3]{merico.argentati[at]cudenver.edu}
\footnotetext[4]{http://math.cudenver.edu/$\sim$aknyazev/, andrew.knyazev[at]cudenver.edu.
Research supported by the NSF-DMS 0612751 and 0728941.}
\footnotetext[5]{School of Computer Science, McGill University, Montreal,
   Quebec, Canada, H3A 2A7, (paige[at]cs.mcgill.ca).
   Research supported by NSERC of Canada Grant OGP0009236.}
\footnotetext[6]{Department of Mathematics and Statistics,
McGill University, Montreal, Quebec, Canada, H3A 2K6,
(ipanay[at]math.mcgill.ca). Research supported by FQRNT of Quebec
Scholarship.}

\renewcommand{\thefootnote}{\arabic{footnote}}

\begin{abstract}
The Rayleigh-Ritz method is widely used for eigenvalue
approximation. Given a matrix $X$ with columns that form an
orthonormal basis for a subspace $\X$, and a Hermitian matrix $A$,
the eigenvalues of $X^HAX$ are called Ritz values of $A$ with
respect to $\X$. If the subspace $\X$ is $A$-invariant then the Ritz
values are some of the eigenvalues of $A$.
If the $A$-invariant subspace $\X$ is perturbed to give
rise to another subspace $\Y$,
then the vector of absolute values of changes in Ritz values of $A$
represents the absolute eigenvalue approximation error using $\Y$.
We bound the error in terms of principal angles between $\X$ and $\Y$.
We capitalize on ideas from a recent paper [DOI: 10.1137/060649070] by
 A. Knyazev and M.~Argentati,
where the vector of absolute values of differences
between Ritz values for subspaces $\X$ and $\Y$ was
weakly (sub-)majorized by a constant times
the sine of the vector of principal angles between $\X$ and $\Y$,
the constant being the spread of the spectrum of $A$.
In that result
no assumption was made on either subspace being $A$-invariant.  It was
conjectured there that if one of the trial
subspaces is $A$-invariant then an analogous weak majorization bound
should be much stronger as it should
only involve terms of the order of sine squared.
Here we confirm this conjecture.  Specifically we prove
that the absolute eigenvalue error is weakly majorized
by a constant times the sine squared of the vector of principal angles
between the subspaces $\X$ and $\Y$, where the constant is
proportional to the spread of the spectrum of $A$.
For many practical cases we show that the proportionality
factor is simply one, and that this bound is sharp. For the general
case we can only prove the result with a slightly larger
constant, which we believe is artificial.
\end{abstract}

\begin{AM}
15A18, 
15A42, 
15A57, 
15A60. 
\end{AM}

(Place for Digital Object Identifier, to get an idea of the final spacing.)  

\pagestyle{myheadings} \thispagestyle{plain}
\markboth{M. E. ARGENTATI, A. V. KNYAZEV, C. C. PAIGE, and I. PANAYOTOV}
{BOUNDS ON CHANGES IN RITZ VALUES}


\section{Introduction} \label{sec:1intro}

Eigenvalue problems appear in many applications.
For example eigenvalues represent the
frequencies of vibration in mechanical vibrations,
while the energy levels of a system are the
eigenvalues of the Hamiltonian
in quantum mechanics.
Eigenvalue problems are used today in these and many other applications,
including spectral data clustering and internet search engines.

Eigenvalues cannot be computed exactly except in some trivial cases,
so numerical approximation is required. Eigenvalue \emph{a
posteriori} and \emph{a priori} error bounds
describe the 
eigenvalue approximation quality, and this is a classical and
important topic in matrix analysis. A posteriori bounds are based on
information readily computable, e.g., the eigenvector residuals, and
are necessary, e.g., for adaptive numerical methods for eigenvalue approximation. 
A priori bounds are given in terms of theoretical
properties, and can be very useful in assessing relative performance
of algorithms.

The widely used Rayleigh-Ritz method is well known for its ability
to generate high quality approximations to eigenvalues of Hermitian
matrices. It is the basis for many numerical procedures for
computing eigenvalues, such as finite element methods and the
Lanczos eigenproblem iteration. Eigenvalue error bounds for the
Rayleigh-Ritz method are important 
since they provide estimates and predictions of the quality of
eigenvalue approximations, and can be used, e.g., to predict the
number of iterations needed in the Lanczos method for computing some
eigenvalues to within a given accuracy. There is a vast literature
on Rayleigh-Ritz eigenvalue methods and error bounds, see, e.g.,\
\cite[Chapter 4]{kvzrs}, \cite[Chapters 10--13]{Par80}, and
\cite[Chapters 3-5]{MR0400004}.

We contribute to this traditional area of research with a new
twist---using weak majorization. Majorization is a classical
technique that can be used
to formulate and prove a great variety of inequalities in a concise
and elegant way. It is widely used in matrix analysis, e.g.,\ to
bound perturbations of eigenvalues via Lidskii's beautiful theorem
\cite{Lid50}. In the context of Rayleigh-Ritz eigenvalue error
bounds, weak majorization was introduced in the celebrated work of Davis and Kahan
\cite{DavKah} to bound eigenvalue errors \emph{a posteriori}. In
the present paper we propose and prove what appear to be the first
theorems based on weak majorization for \emph{a priori}
Rayleigh-Ritz eigenvalue error bounds.
Our results provide a theoretical foundation that can be applied in 
a number of situations, e.g.,\ for finite element methods, e.g.,\ \cite{d}, and for
block Lanczos iterations such as in \cite{GolU77}, see \cite{ka07}. 

We use several well known majorization results found, e.g.,\ in
\cite{Bha97,Horn2,MarO79}. We give references throughout the paper
for the concepts we introduce. For a more thorough background and
reference list, see \cite{KnyA06}.

The rest of the paper is organized as follows.
Section \ref{sec:Prereq} contains all necessary definitions
and basic facts on majorization
that we need for our eigenvalue and singular value bounds.
Section \ref{sec:intro} is the main part of the paper,
where we motivate and formulate our conjectures and theorems.
Section \ref{sec:pertbds} has all our proofs. In section \ref{sec:diss}
we show that our main results are sharp; we also discuss our proofs,
and the possibility that our bound for the most general case might be slightly improved.


\section{Definitions and Prerequisites} \label{sec:Prereq}

We introduce the definitions and tools we need,
together with some mild motivation. We do not provide proofs for the
results in this section---instead we refer the reader to some of the
relevant literature.

\subsection{Notation}\label{sec:notation}

For a real vector $x=[x_1,\ldots, x_n]^T$, we use
$x^\da\equiv[x_1^\da,\ldots, x_n^\da]^T$ to denote $x$ with its
elements rearranged in descending order, while
$x^\ua\equiv[x_1^\ua,\ldots, x_n^\ua]^T$ denotes $x$ with its
elements rearranged in ascending order. We use $|x|$ to denote the
vector $x$ with the absolute value of its components. We use the
`$\leq$' symbol to compare real vectors component-wise.
For real vectors $x$ and $y$ the expression
$x\prec y$ means that $x$ is majorized by $y$,
while $x\wsm y$ means that $x$ is weakly (sub-)majorized by
$y$, see section~\ref{sec:maj}.

We consider the Euclidean space $\C^n$ of column vectors equipped
with the standard scalar product $x^Hy$ and the norm
$\|x\|=\sqrt{x^Hx}$.
We use the same notation $\|A\|$ for the induced matrix
norm 
of a complex matrix $A \in \C^{n\times n}$.
$\X=\Ra(X)\subset \C^n$ means the subspace $\X$ is equal to the
range of the matrix $X$ with $n$ rows. The unit matrix is $I$ and
the zero matrix (not necessarily square) is $0$,
while $e=[1,\ldots,1]^T$. We use $\hm$ to
denote the set of $n\times n$ Hermitian matrices and $\um$ to denote
the set of $n\times n$ unitary matrices in the set $\C^{n\times n}$
of all $n\times n$ complex matrices.

We write $\lambda(A)=\lambda^\da(A)$ for the vector of eigenvalues
of $A \in \hm$ arranged in descending order, and we write $s(B)=
s^\da(B)$ for the vector of singular values of $B$ arranged in
descending order. Individual eigenvalues and singular values are
denoted by $\lambda_i(A)$ and $s_i(B)$, respectively, so, e.g.,\
$\spr(A)=\lambda_1(A)-\lambda_n(A)$ and $s_1(B)=\|B\|$.

Let subspaces $\X$ and $\Y\subseteq \C^n$ have the same dimension, with orthonormal
bases given by the columns of the matrices $X$ and $Y$ respectively.
We denote the vector
of principal angles between $\X$ and $\Y$ arranged in descending order by
$\theta(\X,\Y)=\theta^\da(\X,\Y)$, and define it using
$\cos \theta (\X,\Y)=s^\ua(X^HY)$, e.g.,\
\cite{BjoG73}, \cite[\S12.4.3]{GolL89}.


\subsection{Majorization and Weak Majorization} \label{sec:maj}

We now briefly define the concepts of majorization and weak
majorization which are comparison relations between two real
vectors. For detailed information we refer the reader
to \cite{Bha97,Horn2,MarO79}.

We say that $x\in\R^n$ is weakly (sub-)majorized by $y\in\R^n$,
written $x\wsm y$, if
\be
\sum_{i=1}^k x_i^\da \leq \sum_{i=1}^k y_i^\da,\qquad 1\leq
k\leq n,
\label{eq:wm}
\ee
while $x$ is (strongly) majorized by $y$, written $x\prec y$,
if (\ref{eq:wm}) holds together with
\be \sum_{i=1}^n x_i = \sum_{i=1}^n y_i.
\label{eq:maj}
\ee

Our final results in the paper are weak majorization bounds of the
form $x\prec_w y$ with $x\geq 0$.  On the one hand, we can see from
(\ref{eq:wm}) that $x\leq y \Rightarrow x\wsm y$, i.e., the
inequality implies weak majorization. In our case the 
advantage of using weak majorization is that the inequality $x\leq
y$ (the values of $x$ and $y$ become apparent later) is simply
wrong, while the weak majorization bound $x\prec_w y$ does hold. On
the other hand, a weak majorization bound $x\prec_w y$ implies
$\max(x)\leq \max(y)$. So if the bound $\max(x)\leq\max(y)$ is
already known, but it is also known that $x\leq y$ does not hold,
it makes sense to conjecture and to try to prove $x\prec_w y$ .

Strong `$\prec$' and weak `$\prec_w$' majorization relations share
only some properties with the usual inequality `$\leq$' relation, so
one should deal with them carefully.
For example `$\prec$' and `$\prec_w$' are reflexive and transitive,
but $x\prec y$ and $y\prec x$ do not imply $x=y$; e.g., \cite[Remark II.1.2]{Bha97}.
Similarly $x\prec y$ does not imply the intuitive $x+z\prec y+z$,
as is seen in the example $x=(0,0,0)$, $y=(2,-1,-1)$, $z=(-2,0,0)$.
So we must be particularly careful of the ordering when we \emph{combine} results.
Thus it can be seen from (\ref{eq:wm}) and (\ref{eq:maj}) that:
$x + u \prec x^\da + u^\da$, e.g., \cite[Corollary II.4.3]{Bha97}, and
\be
\{x \prec_w y\}\   \&\   \{u \prec_w v\}\  \&
\ \cdots \ \ \Rightarrow \ \ x + u + \cdots \ \prec \
x^\da + u^\da + \cdots \ \prec_w \
y^\da + v^\da + \cdots,
\label{eq:gen}
\ee
where this also holds with `$\prec_w$' replaced by `$\prec$'.

Some of the other basic majorization and related results we use are fairly obvious:
\begin{align}
&A\in \hm \Rightarrow |\lambda(\pm A)|^\da = s(A);
\label{eq:evsv}\\
&|x\pm y|\wsm |x|^\da + |y|^\da,\quad\mbox{since from (\ref{eq:gen})}\quad
|x\pm y|\leq  |x| + |y| \prec |x|^\da + |y|^\da;
\label{eq:absd}\\
&x\prec y \Rightarrow  |x|\wsm |y|,
\quad\mbox{see, e.g., \cite[Example II.3.5]{Bha97}.}
\label{eq:abs}
\end{align}

Arithmetic operations, e.g.,\ the sum and the product,
on vectors used in majorization are performed component-wise.
In the subsequent Theorems \ref{thm:product} and \ref{thm:norm}
for rectangular matrices we may need to operate with nonnegative vectors of
different lengths. A standard agreement in this case is to add zeroes
at the end of the shorter vector to match the sizes needed for
component-wise arithmetic operations and comparisons.
We also use this agreement in later proofs.

Many inequality relations between eigenvalues and singular values
are succinctly expressed as majorization or weak majorization relations;
and a beautiful example is
%
\begin{theorem}\label{thm:Lid}
\textup{(Lidskii \cite{Lid50}, see also, e.g.,\ \cite[p. 69]{Bha97})}. Let $A$ and $B \in \hm$. The
eigenvalues of $A$, $B$, and $A-B$ satisfy
$\lambda(A)-\lambda(B)\prec \lambda(A-B)$.
\end{theorem}

Recall here that $\lambda(A)-\lambda(B)=\lambda^\da(A)-\lambda^\da(B)$.
Note that the equivalent of (\ref{eq:maj}) holds here using ${\rm trace}(A)=\sum_i \lambda_i(A)$. %
We will use the following corollary:
\begin{corollary} \label{cor:sum}
\textup{(E.g.,\  \cite[Chapter 9, G.1.d]{MarO79},
\cite[Corollary 3.4.3]{Horn2}).} If $A$ and $B \in \C^{n\times n}$ then $s(A \pm B)\wsm
s(A)+s(B)$.
\end{corollary}
This corollary also follows from a weaker statement than
Lidskii's theorem, e.g.,\
\cite[Exercises II.1.14, II.1.15]{Bha97}.

By using (\ref{eq:gen}) we can see that
Corollary \ref{cor:sum} extends to the case of three or more matrices,
because all vectors $s(A), s(B), \ldots$ are nonincreasing.

We also use results for the singular values of a product of matrices:
\begin{theorem}\label{thm:product}
\textup{(E.g.,\  \cite[Theorem 3.3.14]{Horn2}).}
$s(AB)\wsm s(A)s(B)$
for arbitrary, possibly rectangular, matrices $A$ and $B$
such that $AB$ exists.
\end{theorem}
\begin{theorem}\label{thm:norm}
\textup{(E.g.,\
\cite[Theorem 3.3.16]{Horn2}, \cite[Problem III.6.2]{Bha97}).}\\
$s(AB)\leq \|A\|s(B)$ and $s(AB)\leq\|B\|s(A)$
for arbitrary, possibly rectangular, matrices $A$ and $B$
such that $AB$ exists.
\end{theorem}

\section{Motivation and Main Results} \label{sec:intro}

The Rayleigh-Ritz method
for approximating eigenvalues
of a Hermitian matrix $A$ finds the eigenvalues of $X^HAX$, where the
columns of the matrix $X$ form an orthonormal basis for a subspace
$\X$. Here $\X$ is called a trial subspace. The eigenvalues of
$X^HAX$ do not depend on the particular choice of basis and are
called Ritz values of $A$ with respect to $\X$. If $\X$ is
one-dimensional and spanned by the unit vector $x$ there is only one
Ritz value---namely the Rayleigh quotient $x^HAx$. 

When the trial subspace $\X$ is perturbed to become the subspace $\Y$, it is
useful to know how the Ritz values of $A$ vary.
For one-dimensional $\X$ and $\Y$, spanned by unit vectors $x$ and
$y$ respectively, the following result appears in, e.g., \cite[Theorem 1]{ka03}:
\be
|x^HAx-y^HAy|\leq \spr(A)\sin {\theta}(x,y).
\label{eq:knya1}
\ee
Here and below
$\theta(x,y)$ is the acute angle between
the two unit vectors $x$ and $y$ defined by $\theta(x,y)=\mathrm{arccos}
{|x^Hy|}\in [0,\pi/2]$.

It is well known that every eigenvector is a stationary point of the
Rayleigh quotient (considered as a function of a vector)---i.e.,
in the vicinity of an eigenvector the Rayleigh quotient changes very slowly.
The classic result that motivates this paper is the following:
the Rayleigh quotient approximates an eigenvalue of a
Hermitian matrix with accuracy proportional to the \emph{square} of
the eigenvector approximation error. The following simple
bound, e.g.,\ \cite[Theorem 4]{ka03}, demonstrates this:
 \be
|x^HAx-y^HAy|\leq
\spr(A)\sin^2 \theta(x,y),
\label{eq:knya2}
\ee
where we assume that one of the unit vectors $x$ or $y$
is an eigenvector of $A$.
To give a thorough background to our results
we re-derive this important basic bound.
Let $Ax=x\lambda$, then $x^HAx=\lambda$ so
$|x^HAx-y^HAy|=|y^H(A-\lambda I)y|$. We now
plug in the orthogonal decomposition
$y=u+v$ where $u\in\mathrm{span}\{x\}$ and $v\in(\mathrm{span}\{x\})^\perp$.
Thus $(A-\lambda I)u=0$ and $\|v\|=\sin\theta(x,y)$, which results in
$|y^H(A-\lambda I)y|=|v^H(A-\lambda I)v|\leq \|A-\lambda I\| \|v\|^2=
\|A-\lambda I\|\sin^2\theta(x,y)$.
But $\|A-\lambda I\|\leq\spr(A)$,
giving (\ref{eq:knya2}).

Let us now discuss some generalizations of (\ref{eq:knya1}) and
(\ref{eq:knya2}) for subspaces $\X$ and $\Y$ of dimensions higher
than one, with $\dim{\X}=\dim{\Y}$. Let $X$ and $Y$ be two matrices
whose columns form orthonormal bases for $\X$ and $\Y$ respectively,
and suppose that the Ritz values of $A$ with respect to $\X$ and
$\Y$ are arranged in descending (more precisely ``nonincreasing'') order.
To generalize (\ref{eq:knya1}) and (\ref{eq:knya2}) we replace the usual notion of
angles between vectors by a more general one of principal angles
between subspaces, and replace the inequality symbol by the weak
(sub-)majorization symbol `$\prec_w$'.

Let $\lambda(A)$ denote the vector of descending eigenvalues $\lambda_i(A)$ of
a Hermitian matrix $A$, $s(B)$ the vector of descending singular values
of a matrix $B$,
and $\theta (\X,\Y)$ the vector of descending  principal angles $\theta_i(\X,\Y)$
between the subspaces $\X$ and $\Y$, defined such that the vectors
$\cos \theta (\X,\Y)$ and $s(X^HY)$ are the same, except for the reversed order,
see, e.g., \cite{BjoG73}, \cite[Section 12.4.3]{GolL89}. 
A recent paper \cite{KnyA06} generalizes (\ref{eq:knya1}) to:
\be
|\lambda(X^HAX)-\lambda(Y^HAY)|\wsm \spr(A)\sin \theta (\X,\Y).
\label{eq:knya}
\ee
The weak majorization bound (\ref{eq:knya}) implies, e.g.,\ a bound for
its largest term:
\be
\max_i |\lambda_i(X^HAX)-\lambda_i(Y^HAY)| \leq
\spr(A)\,\gap(\X,\Y),
\label{eq:knyamax}
\ee
where
$\gap(\X,\Y)= \max_i \{\sin \theta_i(\X,\Y)\}$
in this case,
e.g.,\ \cite{ka02,KnyA06}.

Both bounds (\ref{eq:knya}) and (\ref{eq:knyamax}) generalize
(\ref{eq:knya1}) to multidimensional subspaces,
but no assumption of $A$-invariance is made in either case. What is the bound that
generalizes (\ref{eq:knya2}), assuming that one of the subspaces
$\X$ or $\Y$ is $A$-invariant? A natural conjecture, made in
\cite{KnyA06}, is that such a bound could be obtained in terms of
$\sin^2\theta (\X,\Y)$. No majorization result of this kind is known,
but simpler results---for the largest error only---are available;
e.g.,\ the following important bound is proved in \cite{k86},
reproduced in \cite[Theorem 2, p. 477]{d}, and \cite[Theorem 2.4]{MR2206452},
with a different proof suggested in \cite[Theorem 2.2.3, p. 56]{k}; for
an English translation of the latter see \cite[Theorem 2.3, p. 383]{ks}.
We present here a slightly modified formulation to make it
consistent with (\ref{eq:knyamax}): if $\X$ or $\Y$ is $A$-invariant
and corresponds to a contiguous set of the extreme, i.e., largest or
smallest, eigenvalues of $A$, then
\be
\max_i |\lambda_i(X^HAX)-\lambda_i(Y^HAY)| \leq \spr(A)\,\gap^2(\X,\Y).
\label{eq:knyamax1}
\ee

Bound (\ref{eq:knyamax1}) generalizes (\ref{eq:knya2}), but does not
take advantage of majorization. Comparing (\ref{eq:knya}) and
(\ref{eq:knyamax}) with (\ref{eq:knya1}), and (\ref{eq:knyamax1})
with (\ref{eq:knya2}), we make an educated guess
for the general case where the invariant 
subspace is not necessarily associated with a contiguous set of extreme eigenvalues:
\begin{conj} \label{conj:main}
Let the subspaces $\X$ and $\Y$ have the same dimension, with orthonormal
bases given by the columns of the matrices $X$ and $Y$ respectively.
Let the matrix $A$ be Hermitian, and $\X$ or $\Y$ be $A$-invariant. Then
\be
|\lambda(X^HAX)-\lambda(Y^HAY)|\wsm
\spr(A)\sin^2 \theta (\X,\Y).
\label{eq:knyamain}
\ee
\end{conj}

We emphasize that the bound (\ref{eq:knyamain}) involves the sine \emph{squared}
and, since convergence analyses are of particular interest for small angles,
this is a great improvement over (\ref{eq:knya}).
This is just as we would hope,
since one of the subspaces is $A$-invariant in (\ref{eq:knyamain}).
The exact $A$-invariance assumption is equivalent to the
subspace being spanned by some exact eigenvectors of $A$,
and Conjecture \ref{conj:main} is an a~priori Rayleigh-Ritz eigenvalue error bound
which can be used to examine how the subspaces $\Y$ of an iterative
eigenproblem algorithm approach an ideal $A$-invariant subspace $\X$.
As we mentioned in the introduction, 
eigenvalue error bounds are important in many applications.
We refer the reader to the follow-up paper \cite{ka07} where
we extend some results of this paper to Hilbert spaces, and discuss
in detail applications to finite element methods and subspace iterations.

The implications of the weak majorization inequality  (\ref{eq:knyamain})
in Conjecture \ref{conj:main} may not be obvious to every reader.
The weak majorization bound (\ref{eq:knyamain}) directly implies
\[
\sum_{i=1}^j |\lambda_i(X^HAX)-\lambda_i(Y^HAY)|^\downarrow
\leq \spr(A) \sum_{i=1}^j \sin^2(\theta_i(\X,\Y))^\downarrow
,\quad \ j=1,\ldots,k,
\]
see (\ref{eq:wm}), where $k=\dim{\X}=\dim{\Y}$.  For example for $j=k$ we obtain
\[
\sum_{i=1}^k |\lambda_i(X^HAX)-\lambda_i(Y^HAY)|
\leq \spr(A) \sum_{i=1}^k \sin^2(\theta_i(\X,\Y)),
\]
and for $j=1$ we get (\ref{eq:knyamax1}). Moreover,
for real vectors $x$ and $y$ the weak majorization $x \prec_w y$
is equivalent to the inequality
$\sum _{i=1}^n \phi(x_i) \leq \sum_{i=1}^n \phi(y_i)$
holding for any continuous nondecreasing convex real valued function $\phi$, see,
e.g.,\ \cite[Statement 4.B.2]{MarO79}.
If for example we take
$\phi(t)=t^p$ with $p \geq 1$,
the bound (\ref{eq:knyamain}) also implies
\[
\left(\sum_{i=1}^k |\lambda_i(X^HAX)-\lambda_i(Y^HAY)|^p\right)^{\frac{1}{p}}
\leq \spr(A) \left(\sum_{i=1}^k \sin^{2p}(\theta_i(\X,\Y))\right)^{\frac{1}{p}}.
\]

We have not proven that Conjecture \ref{conj:main}
holds in \emph{all} circumstances,
and indeed it might not (but we suspect it does).
But we have proven it \emph{always} holds if we multiply the bound by $1.5$.
In section \ref{sec:pertbds} we also show that Conjecture  \ref{conj:main}
does hold in some very useful circumstances:
\begin{theorem} \label{thm:mainex}
The bound (\ref{eq:knyamain}) of Conjecture \ref{conj:main} holds if,
in addition to the assumptions of Conjecture \ref{conj:main},
either or both of the following conditions hold:
\begin{description}
\item[(a)] The $A$-invariant subspace $\X$ or $\Y$ corresponds to a contiguous
set of the largest (or smallest) eigenvalues of $A$.
\item[(b)] All the eigenvalues of $A$ corresponding to the $A$-invariant subspace $\X$ or $\Y$
lie between (and possibly including) one extreme eigenvalue of $A$ and the midpoint
$[\lambda_1(A)\!+\!\lambda_n(A)]/2$ of $A$'s spectrum.
\end{description}
\end{theorem}
This does not cover all known cases where (\ref{eq:knyamain}) holds,
but it does cover many practical cases.
For example in approximating the eigenvalues of a Hermitian matrix,
perhaps using Lanczos' eigenvalue algorithm, e.g., \cite[\S 9]{GolL89},
we are often interested in just one end of the spectrum.
In section \ref{sec:pertbds} we also show a weaker result \emph{always} holds:
\begin{theorem} \label{thm:mainin}
Under the assumptions of Conjecture \ref{conj:main}
we have
\be
|\lambda(X^HAX)-\lambda(Y^HAY)|  \wsm
\spr(A) \left[\ones -\cos \theta (\X,\Y)+\frac{1}{2}
\sin^2\theta (\X,\Y)\right].
\label{eq:bound}
\ee
\end{theorem}
Here and below we use `$\ones$' to indicate a vector of ones.
Note that the individual elements for both vectors
$\ones-\cos\theta(\X,\Y)$ and $\sin^2\theta(\X,\Y)$ are decreasing,
since both functions $1 \!-\!\cos \theta$ and $\sin^2\theta$
are monotonically increasing within $[0,\pi/2]$,
and the vector $\theta(\X,\Y)$ is chosen to be decreasing. %
We now deduce two simple corollaries of Theorem~\ref{thm:mainin}.
Using elementary trigonometry, for $\theta \in [0,\pi/2]$:
\begin{align*}
2-2\cos{\theta}&=2-2\cos{\theta}-(1-\cos{\theta})^2+
(1-\cos{\theta})^2 \nonumber\\
&=\sin^2{\theta}+(1-\cos{\theta})^2
= \sin^2{\theta}+{\sin^4{\theta}}/{(1+\cos{\theta})^2} \nonumber\\
&\leq  \sin^2{\theta}+\sin^4{\theta}. \nonumber
\end{align*}
We first conclude that bound (\ref{eq:bound})
is slightly worse than bound (\ref{eq:knyamain}) from Conjecture
\ref{conj:main}; and second, we immediately obtain from (\ref{eq:bound}):
\begin{corollary}
Under the assumptions of Conjecture \ref{conj:main},
we have
\begin{align}
|\lambda(X^HAX)-\lambda(Y^HAY)| & \wsm
\spr(A) \left[\sin^2 \theta (\X,\Y)+\frac{1}{2}
\sin^4\theta (\X,\Y)\right]
\label{eq:bound1a} 
\\
& \leq \ \ \ \frac{3}{2} \spr(A) \sin^2 \theta (\X,\Y).
\label{eq:bound1}
\end{align}
\end{corollary}
Extending the above trigonometric relation we see that
$$
2-2\cos{\theta}
=\sin^2\theta\left(1+\frac{\sin^2\theta}{(1+\cos\theta)^2}\right)
=\frac{2\sin^2\theta}{1+\cos\theta}
=\frac{\sin^2{\theta}}{\cos^2({\theta/2})}
\leq \tan^2{\theta},
$$
for $\theta \in [0,\pi/2]$; and with $\sin^2{\theta}\leq \tan^2{\theta}$,
bound (\ref{eq:bound}) implies another corollary:
\begin{corollary}
Under the assumptions of Conjecture \ref{conj:main},
we have
\be
|\lambda(X^HAX)-\lambda(Y^HAY)|  \wsm \spr(A)
\tan^2\theta (\X,\Y).
\label{eq:boundtan}
\ee
\end{corollary}

We give an example in section \ref{sec:diss} demonstrating that
the conjectured bound (\ref{eq:knyamain}) cannot be any tighter.
Our numerical tests suggest that
Conjecture \ref{conj:main} holds, i.e., that bound (\ref{eq:bound})
can probably be improved to (\ref{eq:knyamain}).
However we show in section \ref{sec:diss} that already the first step in our proof
of  Theorem \ref{thm:mainin} does not allow us to prove the better
bound (\ref{eq:knyamain}), so a completely different approach is
apparently needed to support Conjecture \ref{conj:main} in all cases---see
section \ref{sec:diss} for more thoughts on this.

Conjecture \ref{conj:main} turns out to be easy to formulate,
but hard to prove in its generality. 
We believe 
that the present publication, which
proves  Conjecture \ref{conj:main} in several practically
interesting particular cases and provides slightly
weaker bounds  (\ref{eq:bound})--(\ref{eq:boundtan}) for the general case, is 
important since it serves as a theoretical foundation
for our future work on applications, e.g.,\ \cite{ka07}.
It is also novel---we know of no other case where majorization is used
for \emph{a priori} Rayleigh-Ritz error bounds.
The only somewhat related result known to us is the pioneering work of
\cite{DavKah}, where majorization is applied to bound
eigenvalue errors \emph{a posteriori}.

\section{Proofs}\label{sec:pertbds}

We have all the tools needed to prove our main results Theorem
\ref{thm:mainex} and Theorem \ref{thm:mainin}. At first both proofs
develop along the same lines; later they split.

By the assumptions in the theorems,
$\X$ and $\Y$ are two subspaces of $\C^n$ of the
same dimension $k$, and are the column ranges of matrices $X$ and
$Y$ with orthonormal columns that
are arbitrary up to unitary transformations of their columns.
Using the singular value decomposition we choose such
a pair of matrices $X$ and $Y$ with orthonormal columns so that
$C\equiv X^HY$ is real, square and diagonal,
with the diagonal entries in increasing order.
Thus by the definition of angles between subspaces,
\be
C=\diag\left(s^\ua(X^HY)\right)=\diag\left(\cos\theta (\X,\Y)\right).
\label{eq:ce}
\ee
We arbitrarily complete $X$ and $Y$ to unitary matrices $[X,X_\perp]$, and $[Y,Y_\perp]\in\um$
and consider the
$2\times2$ partition of their unitary product
$[X,X_\perp]^H[Y,Y_\perp]$. By construction of $X$ and $Y$, its
$k\times k$ upper left block is $C$. We denote its
$(n\!-\!k)\times k$ lower left block by $S\equiv(X_\perp)^H Y$.
Since
$[X,X_\perp]^H[Y,Y_\perp]$ is unitary, the entries $C$ and $S$ of its first
block column satisfy $C^2+S^HS=I$.
So $\lambda(S^H S)=\lambda(I\!-\!C^2)=
\ones \!-\!\cos^2 \theta(\X,\Y) = \sin^2 \theta(\X,\Y),$ where
$\ones$ is the vector of ones, and so the
vectors of singular values $s(C)$ and $s(S)$
are closely connected and we derive from this that
\be
\sin\theta(\X,\Y)=[s(S),0,\ldots,0],
\label{eq:s}
\ee
where $\max\{2k\!-\!n,0\}$ zeroes are added on the right-hand side
to match the number $k$ of angles in the vector $\theta(\X,\Y)$
with the number $\min\{k,n\!-\!k\}$ of singular values in the vector
$s(S)$.

Both theorems assume that either $\X$ or $\Y$ is $A$-invariant,
so without loss of generality let $\X$ be $A$-invariant.
Then since $[X,X_\perp]$ is unitary:
$$
[X,X_\perp]^HA\ [X,X_\perp]=\diag(A_{11},A_{22}),
\mbox{ and }
A=[X,X_\perp]\diag(A_{11},A_{22})[X,X_\perp]^H.
$$
Here $X^HAX=A_{11}\in\mathcal{H}(k)$ and
$(X_\perp)^HAX_\perp=A_{22}\in{\mathcal{H}(n-k)}$.
We can now use $Y^H[X,X_\perp]=[C^H,S^H]=[C,S^H]$
to show that
\begin{align}
Y^HAY&=
Y^H \left([X,X_\perp]\diag(A_{11},A_{22})[X,X_\perp]^H\right)Y
= C A_{11} C + S^H A_{22} S.
\label{eq:yay}
\end{align}

The expression we want to bound in Theorems \ref{thm:mainex} and \ref{thm:mainin}
now takes the form
\begin{align}
\lambda(X^HAX)-\lambda(Y^HAY)
&=\lambda(A_{11})-\lambda(CA_{11}C+S^HA_{22}S)
\nonumber\\
=\lambda(A_{11})-&\lambda(CA_{11}C)+\lambda(CA_{11}C)-\lambda(CA_{11}C+S^HA_{22}S)
\nonumber\\
\prec [\lambda(A_{11})-&\lambda(CA_{11}C)]^\da+\lambda(-S^HA_{22}S),
\label{eq:ident}
\end{align}
where this last line used Lidskii's Theorem \ref{thm:Lid} with (\ref{eq:gen}).
See the discussion following (\ref{eq:ex}) for more about this choice.
Next (\ref{eq:evsv}), Theorems \ref{thm:product} and \ref{thm:norm},
and (\ref{eq:s}) give
\be
\left|\lambda(-S^HA_{22}S)\right|^\da = 
s(S^HA_{22}S) \wsm \|A_{22}\|\sin^2\theta (\X,\Y).
\label{eq:bd22}
\ee

At this point the proofs split.
Each proof will use a different majorization of
$\lambda(A_{11})-\lambda(CA_{11}C)$ in (\ref{eq:ident}),
but both will use (\ref{eq:bd22}).
We first establish Theorem \ref{thm:mainex}.
Neither (\ref{eq:knyamain}) nor (\ref{eq:bound})
is altered by replacing $A$ by $\pm A\!+\!\alpha I$
where $\alpha$ is an arbitrary real constant,
and so we can make the new $A_{11}$ nonnegative definite
in each of the parts \textbf{(a)} and \textbf{(b)} of Theorem~\ref{thm:mainex}
by choosing the appropriate sign and the shift $\alpha$.

\begin{proof}[of Theorem \ref{thm:mainex}]
The starting point of the proof is (\ref{eq:ident}),
but now we assume $A_{11}$ is nonnegative definite
and so has a nonnegative definite square root $\sqrt{A_{11}}$.
We deal with $\lambda(A_{11})-\lambda(CA_{11}C)$ first. For arbitrary square
matrices $F$ and $G$  we have $\lambda(FG)=\lambda(GF)$.
Taking $F=C\sqrt{A_{11}}$ and $G= \sqrt{A_{11}} C$, we get
$\lambda(CA_{11}C)=\lambda(\sqrt{A_{11}}C^2\sqrt{A_{11}})$.
Using this and Lidskii's Theorem \ref{thm:Lid} we see that
\begin{align*}
\lambda(A_{11})-\lambda(CA_{11}C)
&=\lambda(A_{11})-\lambda\left(\sqrt{A_{11}}C^2\sqrt{A_{11}}\right)\\
&\prec \lambda\left(\sqrt{A_{11}}\sqrt{A_{11}}-\sqrt{A_{11}}C^2\sqrt{A_{11}}\right)\nonumber\\
&= \lambda\left(\sqrt{A_{11}}\left(I-C^2\right)\sqrt{A_{11}}\right)
= \lambda\left(\sqrt{A_{11}}S^HS\sqrt{A_{11}}\right),
\end{align*}
since $C^2+S^HS=I$.
Then using (\ref{eq:abs}) with Theorem 
\ref{thm:norm} (twice), and (\ref{eq:s}) 
we obtain
$$
\left| \lambda(A_{11})-\lambda(CA_{11}C) \right|
\wsm s\left(\sqrt{A_{11}}S^HS\sqrt{A_{11}}\right) 
\leq \|A_{11}\|\sin^2\theta(\X,\Y).
$$
Apply (\ref{eq:abs}) to (\ref{eq:ident}); then (\ref{eq:absd}), (\ref{eq:gen}),
and (\ref{eq:bd22}) with the above bound give
\begin{align}
\left|\lambda(X^HAX)-\lambda(Y^HAY)\right|
&\wsm \left| [\lambda(A_{11})-\lambda(CA_{11}C)]^\da+ \lambda(-S^HA_{22}S)\right|
\nonumber\\
&\wsm \left| \lambda(A_{11})-\lambda(CA_{11}C)\right|^\da+
\left| \lambda(-S^HA_{22}S)\right|^\da
\nonumber\\
&\wsm (\|A_{11}\|+\|A_{22}\|)\sin^2\theta(\X,\Y).
\label{eq:normsum}
\end{align}

Here this proof splits,
and we first prove part \textbf{(a)} of Theorem~\ref{thm:mainex}.
By assumption the invariant subspace $\X$ corresponds
to a contiguous set of the largest (or smallest) eigenvalues of $A$.
Here we present the proof for the case of the largest eigenvalues.
The case of the smallest eigenvalues follows immediately by
substituting $-A$ for $A$.
We replace $A$ with $A+\alpha I$ where $\alpha$ is
chosen as the constant real shift that makes the new $A_{11}$
positive semidefinite (nonnegative definite and singular), 
so that $\sqrt{A_{11}}$ exists. Since
$\dim \X = k$, and the invariant subspace  $\X$ corresponds to a
contiguous set of the largest eigenvalues of $A$,
$\alpha= -\lambda_k(A)$. 
After the shift $\lambda_k(A)$ becomes zero, the
eigenvalues of the block $A_{11}$ become nonnegative with
$\lambda_1(A)$ being the largest in absolute value, while the
eigenvalues of the block $A_{22}$ become nonpositive with
$
\|A_{22}\| =  -\lambda_n(A).
$
Thus
$
\|A_{11}\|+\|A_{22}\| = \lambda_1(A) -\lambda_n(A) = \spr(A).
$
Using this together with (\ref{eq:normsum}) gives (\ref{eq:knyamain}),
completing the proof of part \textbf{(a)}.

For part \textbf{(b)} of Theorem~\ref{thm:mainex}
we prove the case where the eigenvalues of $A_{11}$
lie in the top half of the spectrum of $A$,
the remaining case is proven by substituting $-A$ for $A$.
Choose the shift so that for the new $A$, $\lambda_1(A)=-\lambda_n(A)$,
ensuring with the assumptions that $A_{11}$ is nonnegative definite and that
$\|A_{11}\|\leq\spr(A)/2$ and $\|A_{22}\|\leq\spr(A)/2$,
so that (\ref{eq:normsum}) again leads to (\ref{eq:knyamain}).
\end{proof}

In fact whenever we can choose the sign and shift in $\pm A+\alpha I$
so that this new $A$ has $A_{11}$ nonnegative definite with
$\|A_{11}\|+\|A_{22}\|\leq \spr(A)$, then (\ref{eq:knyamain}) will be satisfied.

We return again to
(\ref{eq:ident}) and (\ref{eq:bd22}) to establish Theorem \ref{thm:mainin}.

\begin{proof}[of Theorem \ref{thm:mainin}]
Applying Lidskii's Theorem \ref{thm:Lid} with (\ref{eq:gen}) to (\ref{eq:ident}) gives
\begin{align}
\lambda(X^HAX)-\lambda(Y^HAY)
&\prec [\lambda(A_{11})-\lambda(CA_{11}C)]^\da+\lambda(-S^HA_{22}S)
\nonumber\\
&\prec \lambda(A_{11}-CA_{11}C)+\lambda(-S^HA_{22}S).
\label{eq:ident2}
\end{align}
In order to bound this we will use the identity
\be
A_{11}-CA_{11}C=(I\!-\!C)A_{11}+CA_{11}(I\!-\!C),
\label{eq:ident3}
\ee
together with the following results obtained using (\ref{eq:ce})
with Theorems \ref{thm:norm} 
and \ref{thm:product}:
\begin{align}
s((I-C)A_{11})&\leq \|A_{11}\|s(I-C) =\|A_{11}\|(\ones-\cos\theta(\X,\Y)),
\label{eq:step3term1}\\
s(CA_{11}(I-C)) &\wsm s(C)s(A_{11}(I-C))\leq s(A_{11}(I-C))
\nonumber\\
&\leq \|A_{11}\|s(I-C)=\|A_{11}\|(\ones-\cos\theta(\X,\Y)).
\label{eq:step3term2}
\end{align}
Discarding the first $C$ in $s(CA_{11}(I-C))$ is no real loss,
see  section~\ref{sec:diss}.
Using (\ref{eq:ident3}) 
and applying (\ref{eq:evsv}),
Corollary \ref{cor:sum}, and (\ref{eq:gen})
with (\ref{eq:step3term1}) and (\ref{eq:step3term2}), gives
\begin{align}
 \left|\lambda(A_{11}-CA_{11}C)\right|^\da
&= s((I-C)A_{11}+CA_{11}(I-C))
\nonumber\\
&\wsm s((I-C)A_{11})+s(CA_{11}(I-C))
\nonumber\\
&\wsm 2\|A_{11}\|(\ones-\cos\theta(\X,\Y)). \label{eq:step3}
\end{align}
Now apply (\ref{eq:abs}) to (\ref{eq:ident2}), followed by (\ref{eq:absd}),
and use (\ref{eq:step3}) and (\ref{eq:bd22}) with (\ref{eq:gen}),
together with $\|A_{11}\|, \|A_{22}\|\leq \|A\|$, to obtain:
\begin{align}
\left|\lambda(X^HAX)-\lambda(Y^HAY)\right|
&\wsm \left| \lambda(A_{11}-CA_{11}C)+\lambda(-S^HA_{22}S)\right|
\nonumber\\
&\wsm  \left|\lambda(A_{11}-CA_{11}C)\right|^\da+\left|\lambda(-S^HA_{22}S)\right|^\da
 \nonumber\\
&\wsm \|A\|\left[2(\ones  - \cos{\theta(\X,\Y)})+
\sin^2{\theta(\X,\Y)}\right].
\label{eq:step4}
\end{align}

Our final step is to replace $\|A\|$ by an expression involving
$\spr(A)$. Observe here that the difference between Ritz values is
invariant under \emph{any} shift $\alpha \in \R$. So we shift $A$ in a way
to minimize $\|A\|$. This situation occurs when 0 is exactly in the
middle of the spectrum, in which case $\|A\|=\spr(A)/2$. Combining
this observation with (\ref{eq:step4}) completes the proof of
(\ref{eq:bound}).
\end{proof}

\section{Discussion}\label{sec:diss}

The following example shows that the conjectured bound (\ref{eq:knyamain})
cannot be improved as a general result.
Let $n=2m$ and let an arbitrary set of $m$ angles $\theta_i$ be given,
where $\pi/2 \geq\theta_1\geq \ldots\geq\theta_{m}\geq 0$.
Let $C=\diag(\cos(\theta_1),\ldots,\cos(\theta_{m}))$,
$X=[I,0]^H$, $Y=\left[C,\sqrt{I-C^2}\right]^H$,
and
$
A=\left[\begin{array}{cc}
I & 0\\
0 & -I\
\end{array}\right],
$
where all unit matrices $I$ are of size $m$, so that $X$ and $Y$
are $n\times m$ and $A$ is $n\times n$. Then the $\theta_i$ become
the principal angles between the pair of $k=m$ dimensional
subspaces $\X\equiv\Ra(X)$ and $\Y\equiv\Ra(Y)$. Moreover the Ritz
values are the eigenvalues of $X^HAX=I$ and $Y^HAY=2C^2-I$, and so
$\left|\lambda(X^HAX)-\lambda(Y^HAY)\right|^\da= 2\sin^2 \theta (\X,\Y)$.
In this example $\spr(A) = 1 - (-1) = 2$, so
(\ref{eq:knyamain}) turns into an equality.

Asymptotically where all of the angles are small, bounds
(\ref{eq:knyamain}), (\ref{eq:bound}), 
(\ref{eq:bound1a}) and (\ref{eq:boundtan}) are all equivalent.
Moreover our numerical tests support Conjecture \ref{conj:main} in
all cases. Perhaps in practical terms, from the point of view of a numerical analyst
we are done. However it would be pleasing to
know whether Conjecture \ref{conj:main} holds theoretically in its
generality, since bound (\ref{eq:knyamain}) looks more aesthetic and
cannot be improved as a general result.

One important thing we know is that our approach of starting with
Theorem \ref{thm:Lid} to deduce (\ref{eq:ident}) (used
in the proof of (\ref{eq:bound})) cannot reduce bound
(\ref{eq:bound}) to bound (\ref{eq:knyamain}) in general, no matter how we
modify the rest of the proof. This can be seen from
the following example in $\C^4$. Let $A=\diag(A_{11},A_{22})$ and
$C=X^HY$, $S=X_\perp^HY$ be as in
\be
A=\left[
\begin{array}{cc|cc}
0&1&0&0 \\
1&0&0&0\\
\hline
0&0&1&0\\
0&0&0&1\\
\end{array}
\right], \ \ [X,X_\perp]=I_4,
\ \ [Y,Y_\perp]=\left[
\begin{array}{cc|cc}
0&0&-1&0 \\
0&1&0&0\\
\hline
1&0&0&0\\
0&0&0&1\\
\end{array}
\right],
\label{eq:ex}
\ee
where $I_4$ is the $4\times 4$ unit matrix, so that
$X,X_\perp,Y,Y_\perp$ all have two columns.
Then $[X,X_\perp]^H[Y,Y_\perp] = [Y,Y_\perp]$ are chosen
as in our proofs
and we see that $\theta (\X,\Y) = [\pi/2,0]^T$,
$CA_{11}C=0$, $S^HA_{22}S=\diag(1,0)$.
Here the
largest and smallest eigenvalues of $A$ are $\pm 1$, so
$\spr(A)=2$.  Hence
by direct calculation
$$
X^HAX = A_{11}=
\left[\ba {cc} 0 & 1 \\
 1 &  0
\ea \right],\quad
Y^HAY =CA_{11}C+S^HA_{22}S=S^HA_{22}S=
\left[\ba {cc} 1 & 0 \\
 0 &  0
\ea \right],
$$
$$
\left|\lambda(X^HAX)-\lambda(Y^HAY)\right| =
\left|\begin{bmatrix} 1\\ -1 \end{bmatrix} - \begin{bmatrix} 1\\ 0 \end{bmatrix}\right| =
\begin{bmatrix} 0\\ 1 \end{bmatrix} \wsm \begin{bmatrix} 2\\ 0 \end{bmatrix}
=\spr(A)\sin^2{\theta (\X,\Y)},
$$
so 
example (\ref{eq:ex}) \emph{does} satisfy (\ref{eq:knyamain}).

Let us now attempt to use (\ref{eq:ident}) for (\ref{eq:ex}). %
The right hand side of (\ref{eq:ident}) is
$$
a\equiv [\lambda(A_{11})-\lambda(CA_{11}C)]^\da+\lambda(-S^HA_{22}S)=
\begin{bmatrix} 1 \\-1\end{bmatrix}-\begin{bmatrix} 0\\ 0\end{bmatrix}+
\begin{bmatrix} 0 \\ -1\end{bmatrix}=\begin{bmatrix} 1 \\-2 \end{bmatrix},
$$
where it is \emph{not} true that $|a|\wsm \spr(A)\sin^2{\theta(\X,\Y)}$.
That is, the absolute value of
the right-hand side of (\ref{eq:ident}) is not always weakly
majorized by $\spr(A)\sin^2{\theta(\X,\Y)}$,
so we cannot obtain a general proof of (\ref{eq:knyamain})
starting from the majorization in (\ref{eq:ident}).

Example (\ref{eq:ex}) can tell us even more.
For \emph{any} matrix $M=M^H$ we have
the following generalization of (\ref{eq:ident2})
($M=CA_{11}C$ in (\ref{eq:ident}) and (\ref{eq:ident2})):
\begin{align*}
\lambda(X^HAX)-\lambda(Y^HAY)
&=\lambda(X^HAX)-\lambda(M)+\lambda(M)-\lambda(Y^HAY)\\
&\prec \lambda(X^HAX-M)+\lambda(M-Y^HAY)\equiv \tilde{a}.
\end{align*}
It might be thought that if, e.g.,\ $X^HAX$ is indefinite,
some such $M$ could be chosen to minimize $\tilde{a}$ and to
prove (\ref{eq:knyamain}).
But in example (\ref{eq:ex}) it can be shown
that there is \emph{no} real symmetric $M$
giving $\tilde{a}$ satisfying the desired bound
$|\tilde{a}|\!\wsm\! \spr(A)\sin^2{\theta (\X,\Y)}$.
In particular $M\!=\!Y^H\!AY$ will not give this bound,
as the reader can check via (\ref{eq:ex}).
That is, using $\lambda(X^HAX)-\lambda(Y^HAY)
\prec \lambda(A_{11}-CA_{11}C-S^HA_{22}S)$
in place of (\ref{eq:ident}) will still not give (\ref{eq:knyamain}) via our approach.

So on the one hand we cannot improve bound (\ref{eq:bound}) to give
(\ref{eq:knyamain}) except possibly by considering a different
approach to our present way of using Lidskii's Theorem \ref{thm:Lid}
or equivalent in the first step, see (\ref{eq:ident}) and (\ref{eq:ident2}).
On the other hand our numerical tests
suggest that the tighter bound (\ref{eq:knyamain}) holds.
Thus if we are to prove (\ref{eq:knyamain}) for widely spread
interior eigenvalues, we appear to need an approach more sophisticated than our
particular application of Lidskii's theorem in the first step.

An essentially equivalent first step
was used in \cite[Theorem 10]{ka03} in an earlier attempt to prove (\ref{eq:knya}),
where it led to an artificial multiplier $\sqrt{2}$ in the right-hand side
of (\ref{eq:knya}). The subsequent 
paper \cite{KnyA06} used an
unusual technique to extend an arbitrary Hermitian operator to an
orthogonal projector in a higher dimensional space, preserving its
Ritz values, to prove (\ref{eq:knya}) as it is stated, without the multiplier
$\sqrt{2}$. Perhaps the same technique might shed light here, and
help us to establish Conjecture \ref{conj:main}, but this currently
remains an open question.


\section*{Conclusions}

We clarify a conjecture of Knyazev and Argentati \cite{KnyA06} on a
bound for the absolute difference between Ritz values of a Hermitian
matrix $A$ for two trial subspaces, one of which is $A$-invariant.
We prove the conjecture for the cases where \textbf{(a)}: the $A$-invariant
subspace corresponds to a contiguous set of the largest (or
smallest) eigenvalues of $A$,
and \textbf{(b)}: the eigenvalues of $A$ corresponding to the $A$-invariant
subspace all lie in the top (or the bottom) half of the spectrum of $A$.
We prove a slightly weaker bound for
general invariant subspaces. We believe that the conjecture holds,
i.e., that this weaker bound can be improved, and this is supported by our
numerical tests, but the proof of the conjecture in its generality
(if it is true) 
may require an unorthodox approach, perhaps one such as that
used in \cite{KnyA06}.
These results are useful in practice, and for example are applicable to the analysis of
routines which use the Rayleigh-Ritz method, such as some Krylov subspace methods.
We refer the reader to the subsequent 
paper \cite{ka07}, where we extend some results of this paper to Hilbert spaces,
and discuss in detail
their application to finite element methods and subspace iterations. \vskip12pt

\def\refname{

\end{document}